\documentclass{amsart}

\title[On the $r-$stability of spacelike hypersurfaces]{On the $r-$stability of spacelike hypersurfaces}

\usepackage{graphicx,pstricks,pst-plot}

\newtheorem{theorem}{Theorem}[section]
\newtheorem{lemma}[theorem]{Lemma}
\newtheorem{proposition}[theorem]{Proposition}
\newtheorem{corollary}[theorem]{Corollary}

\theoremstyle{definition}
\newtheorem{definition}[theorem]{Definition}

\theoremstyle{remark}

\numberwithin{equation}{section}

\author{F. Camargo}
\address{Departamento de Matem\'atica e Estat\'{\i}stica,
Universidade Federal de Campina Grande, Campina Grande,
Para\'{\i}ba, Brazil. 58109-970} \email{fernandaecc@dme.ufcg.edu.br}

\author{A. Caminha}
\address{Departamento de Matem\'atica, Universidade Federal do Cear\'a, Fortaleza,
Cear\'a, Brazil. 60455-760} \email{antonio.caminha@gmail.com}

\author{M. da Silva}
\address{Centro de Matem\'{a}tica, Computa\c{c}\~{a}o e Cogni\c{c}\~{a}o, Universidade Federal do ABC, Santo Andr\'e, S\~{a}o Paulo,
Brazil. 09090-400 } \email{marcfab@gmail.com}

\author{H. de Lima}
\address{Departamento de Matem\'atica e Estat\'{\i}stica,
Universidade Federal de Campina Grande, Campina Grande,
Para\'{\i}ba, Brazil. 58109-970} \email{henrique@dme.ufcg.edu.br}



\thanks{The second author is partially supported by CNPq, Brazil. The third
author is partially supported by CNPq/FAPESQ/PPP, Brazil}

\begin{document}

\maketitle

\begin{abstract}
In this paper we study the strong stability of spacelike hypersurfaces with constant $r$-th mean curvature in Generalized Robertson-Walker spacetimes of constant sectional curvature. In particular, we treat the case in which the ambient spacetime is the de Sitter space.
\end{abstract}

\section{Introduction}
The notion of stability concerning hypersurfaces of constant mean
curvature were first studied by Barbosa and do Carmo in
\cite{BdC:84}, and Barbosa, do Carmo and Eschenburg in
\cite{BdCE:88}, where they proved that spheres are the only stable
critical points of the area functional for volume-preserving
variations.

Related to the Lorentz context, in 1993,  Barbosa and Oliker
\cite{Barbosa:93} obtained an analogous result, proving that
constant mean curvature spacelike hypersurfaces in Lorentz manifolds
are also critical points of the area functional for variations
keeping the volume constant. In this sense, the variational methods
for Riemannian and Lorentz manifolds coincide. They computed the
second variation formula and obtained in the de Sitter space
$\mathbb S_1^{n+1}$ that spheres maximize the area functional for
variation keeping the volume constant; this fact determines the
definition of stability which is given below in a more general form.
Then, in \cite{Barbosa:93} is proved that if $M^n$ is a complete
spacelike immersed hypersurface in $\mathbb S_1^{n+1}$ with constant
mean curvature $H,$ then $M^n$ is stable if $M^n$ is compact, or
$H^2\geq 1$, or $H^2<4(n-1)/n^2$. This result was extended for the
case of complete spacelike hypersurfaces with constant $r$-th mean
curvature in de Sitter space by Brasil and Colares, in
\cite{Brasil:03}.

More recently, the second author joint with Barros and Brasil
\cite{BBC:08} have studied the problem referring to strong stability
(that is, stability without the hypothesis of preserving-volume
variations) for spacelike hypersurfaces with constant mean curvature
in a Generalized Robertson-Walker (GRW) spacetime, giving a
characterization of maximal and spacelike slices. In \cite{Liu:08},
Liu and Yang obtained a extension of the result of \cite{BBC:08} for
spacelike hypersurfaces with constant scalar curvature.

Here, motivated by the works \cite{BBC:08} and \cite{Liu:08}, we
consider spacelike hypersurfaces with constant $r$-th mean curvature
in GRW spacetimes of constant sectional curvature in order to
classify the strongly $r$-stable ones. For this, we will use a
formula due to Barros and Sousa \cite{Barros:09} for a operator
$L_r$, naturally attached to the operators $P_r$ that can be defined
using the $r$-th mean curvatures, for a suitable support function.
More precisely, we will prove the following result:

\begin{theorem}
Let $\overline M^{n+1}_c=I\times_{\phi}F^n$ be a generalized
Robertson-Walker spacetime of constant sectional curvature $c$, and
$x:M^n\rightarrow\overline M^{n+1}_c$ be a closed strongly
$r$-stable spacelike hypersurface. If the warping function $\phi$ is
such that $H_r\phi''\geq\max\{H_{r+1}\phi',0\}$ and the set where
$\phi'=0$ has empty interior on $M^n$, then either $M^n$ is $r$-
maximal or a totally umbilical slice $\{s_0\}\times F$.
\end{theorem}

An application of this previous theorem is obtained on de Sitter
space.

\begin{corollary}
Let $x:M^n\rightarrow\mathbb S_1^{n+1}$ be a closed strongly
$r-$stable spacelike hypersurface. Suppose that the set of points in
which $M^n$ intersects the equator of $\mathbb S_1^{n+1}$ has empty
interior. If $H_r\geq\max\{H_{r+1},0\}$, then either $M^n$ is
$r-$maximal or a umbilical round sphere.
\end{corollary}

\section{$r-$stability of spacelike hypersurfaces}

In what follows, $\overline M^{n+1}$ denotes a time-oriented Lorentz
manifold with Lorentz metric $\overline g=\langle\,\,,\,\,\rangle$,
volume element $d\overline M$ and semi-Riemannian connection
$\overline\nabla$. In this context, we consider spacelike
hypersurfaces $x:M^n\rightarrow\overline M^{n+1}$, namely, isometric
immersions from a connected, $n-$dimensional orientable Riemannian
manifold $M^n$ into $\overline M$. We let $\nabla$ denote the
Levi-Civita connection of $M^n$.

If $\overline M$ is time-orientable and $x:M^n\rightarrow\overline
M^{n+1}$ is a spacelike  hypersurface, then $M^n$ is orientable
(cf.~\cite{ONeill:83}) and one can choose a globally defined unit
normal vector field $N$ on $M^n$ having the same time-orientation of
$\overline M$. One says that such an $N$ {\em points to the future}.

In this setting, let $A$ denotes the corresponding shape operator.
At each $p\in M^n$, $A$ restricts to a self-adjoint linear map
$A_p:T_pM\rightarrow T_pM$. For $1\leq r\leq n$, let $S_r(p)$ denote
the $r-$th elementary symmetric function on the eigenvalues of
$A_p$; this way one gets $n$ smooth functions
$S_r:M^n\rightarrow\mathbb R$, such that
$$\det(tI-A)=\sum_{k=0}^n(-1)^kS_kt^{n-k},$$
where $S_0=1$ by definition. If $p\in M^n$ and $\{e_k\}$ is a basis
of $T_pM$ formed by eigenvectors of $A_p$, with corresponding
eigenvalues $\{\lambda_k\}$, one immediately sees that
$$S_r=\sigma_r(\lambda_1,\ldots,\lambda_n),$$
where $\sigma_r\in\mathbb R[X_1,\ldots,X_n]$ is the $r-$th
elementary symmetric polynomial on the indeterminates
$X_1,\ldots,X_n$.

For $1\leq r\leq n$, one defines the $r-th$ mean curvature $H_r$ of
$x$ by
$${n\choose r}H_r=(-1)^rS_r=\sigma_r(-\lambda_1,\ldots,-\lambda_n).$$

A spacelike hypersurface $x:M^n\rightarrow\overline{M}^{n+1}$ such
that $H_{r+1}=0$ on $M^n$ is said to be {\em $r-$maximal}.

For $0\leq r\leq n$ one defines the $r-$th Newton transformation
$P_r$ on $M^n$ by setting $P_0=I$ (the identity operator) and, for
$1\leq r\leq n$, via the recurrence relation
\begin{equation}\label{eq:Newton operators}
P_r=(-1)^rS_rI+AP_{r-1}.
\end{equation}
A trivial induction shows that
$$P_r=(-1)^r(S_rI-S_{r-1}A+S_{r-2}A^2-\cdots+(-1)^rA^r),$$
so that Cayley-Hamilton theorem gives $P_n=0$. Moreover, since $P_r$
is a polynomial in $A$ for every $r$, it is also self-adjoint and
commutes with $A$. Therefore, all bases of $T_pM$ diagonalizing $A$
at $p\in M^n$ also diagonalize all of the $P_r$ at $p$. Let
$\{e_k\}$ be such a basis. Denoting by $A_i$ the restriction of $A$
to $\langle e_i\rangle^{\bot}\subset T_p\Sigma$, it is easy to see
that
$$\det(tI-A_i)=\sum_{k=0}^{n-1}(-1)^kS_k(A_i)t^{n-1-k},$$
where
$$S_k(A_i)=\sum_{\stackrel{1\leq j_1<\ldots<j_k\leq n}{j_1,\ldots,j_k\neq i}}\lambda_{j_1}\cdots\lambda_{j_k}.$$

With the above notations, it is also immediate to check that
$P_re_i=(-1)^rS_r(A_i)e_i$, and hence (Lemma 2.1
of~\cite{Barbosa:97})
\begin{enumerate}
\item[(a)] $S_r(A_i)=S_r-\lambda_iS_{r-1}(A_i)$;
\item[(b)] ${\rm tr}(P_r)=(-1)^r\sum_{i=1}^nS_r(A_i)=(-1)^r(n-r)S_r=b_rH_r$;
\item[(c)] ${\rm tr}(AP_r)=(-1)^r\sum_{i=1}^n\lambda_iS_r(A_i)=(-1)^r(r+1)S_{r+1}=-b_rH_{r+1}$;
\item[(d)] ${\rm tr}(A^2P_r)=(-1)^r\sum_{i=1}^n\lambda_i^2S_r(A_i)=(-1)^r(S_1S_{r+1}-(r+2)S_{r+2})$,
\end{enumerate}
where $b_r=(n-r){n\choose r}$.

Associated to each Newton transformation $P_r$ one has the second
order linear differential operator $L_r:\mathcal
D(M)\rightarrow\mathcal D(M)$, given by
$$L_r(f)={\rm tr}(P_r\,\text{Hess}\,f).$$

According to \cite{Alias:07}, if $\overline M^{n+1}$ is of constant
sectional curvature, then $P_r$ is a divergence-free and,
consequently,
$$L_r(f)={\rm div}(P_r\nabla f).$$

For future use, we recall Lemma 2.6 of~\cite{Caminha:06}: if
$(a_{ij})$ denotes the matrix of $A$ with respect to a certain
orthonormal basis $\beta=\{e_k\}$ of $T_pM$, then the matrix
$(a_{ij}^r)$ of $P_r$ with respect to the same basis is given by
\begin{equation}\label{eq:fundamental_relation_Reilly}
a_{ij}^r=\frac{(-1)^r}{r!}\sum_{i_k,j_k=1}^n\epsilon_{i_1\ldots
i_ri}^{j_1\ldots j_rj} a_{j_1i_1}\ldots a
_{j_ri_r},
\end{equation}
where
$$\epsilon_{i_1\ldots i_r}^{j_1\ldots j_r}=\left\{\begin{array}{lll}
{\rm sgn}(\sigma)&,&\text{if the}\,\,i_k\,\,\text{are pairwise distinct and}\\
&&\sigma=(j_k)\,\,\text{is a permutation of them};\\
0&,&\text{else}.
\end{array}\right.$$

If $x$ is as above, a {\em variation} of it is a smooth mapping
$$X:M^n\times(-\epsilon,\epsilon)\rightarrow\overline M^{n+1}$$
satisfying the following conditions:
\begin{enumerate}
\item[(1)] For $t\in(-\epsilon,\epsilon)$, the map $X_t:M^n\rightarrow\overline M^{n+1}$
given by $X_t(p)=X(t,p)$ is a spacelike immersion such that $X_0=x$.
\item[(2)] $X_t\big|_{\partial M}=x\big|_{\partial M}$, for all $t\in(-\epsilon,\epsilon)$.
\end{enumerate}

In all that follows, we let $dM_t$ denote the volume element of the metric induced on $M$ by $X_t$ and $N_t$ the unit normal vector field along $X_t$.

The {\em variational field} associated to the variation $X$ is the
vector field $\frac{\partial X}{\partial t}\Big|_{t=0}$. Letting
$f=-\langle\frac{\partial X}{\partial t},N_t\rangle$, we get
\begin{equation}\label{eq:decomposition of variational vector field}\frac{\partial X}{\partial t}=fN_t+\left(\frac{\partial X}{\partial
t}\right)^{\top},\end{equation} where $\top$ stands for tangential
components.

The {\em balance of volume} of the variation $X$ is the function $\mathcal V:(-\epsilon,\epsilon)\rightarrow\mathbb R$ given by
$$\mathcal V(t)=\int_{M\times[0,t]}X^*(d\overline M),$$
and we say $X$ is {\em volume-preserving} if $\mathcal V$ is constant. The following lemma is classical (cf.~\cite{Xin:03}).

\begin{lemma}\label{lemma:first variation}
Let $\overline M^{n+1}$ be a time-oriented Lorentz manifold and $x:M^n\rightarrow\overline M^{n+1}$ a closed spacelike hypersurface. If 
$X:M^n\times(-\epsilon,\epsilon)\rightarrow\overline M^{n+1}$ is a variation of $x$, then
$$\frac{d\mathcal V}{dt}=\int_MfdM_t.$$
In particular, $X$ is volume-preserving if and only if $\int_MfdM_t=0$ for all $t$.
\end{lemma}

We remark that Lemma 2.2 of~\cite{BdCE:88} remains valid in the Lorentz context, i.e., if $f_0:M\rightarrow\mathbb R$ is a smooth function such that 
$\int_Mf_0dM=0$, then there exists a volume-preserving variation of $M$ whose variational field is $f_0N$. Moreover, if we drop the requirement that  variation be volume-preserving (or, which is the same, that $\int_Mf_0dM=0$), the argument of that Lemma always gives a variation whose variational field is $f_0N$.

In order to extend~\cite{Barbosa:97} to the Lorentz setting, we define the {\em $r-$area functional} $\mathcal A_r:(-\epsilon,\epsilon)\rightarrow\mathbb R$ associated to the variation $X$ be given by
$$\mathcal A_r(t)=\int_MF_r(S_1,S_2,\ldots,S_r)dM_t,$$
where $S_r=S_r(t)$ and $F_r$ is recursively defined by setting $F_0=1$, $F_1=-S_1$ and, for $2\leq r\leq n-1$,
$$F_r=(-1)^rS_r-\frac{c(n-r+1)}{r-1}F_{r-2}.$$

The next step is the Lorentz analogue of Proposition 4.1 of~\cite{Barbosa:97}. Since it seems to us that their proof only works on a neighborhood free of umbilics, and in order to keep this work self-contained, we present an alternative one here.

\begin{lemma}\label{lemma:computing the time-derivative of Sr}
Let $x:M^n\rightarrow\overline M^{n+1}_c$ be a closed spacelike hypersurface
of the time-oriented Lorentz manifold $\overline M^{n+1}_c$ with constant sectional curvature $c$, and let
$X:M^n\times(-\epsilon,\epsilon)\rightarrow\overline M^{n+1}_c$ be a variation of $x$. Then,
\begin{equation}\label{eq:differentiation of $S_r$}
\frac{\partial S_{r+1}}{\partial t}=(-1)^{r+1}\left[L_{r}f+c{\rm
tr}(P_{r})f-{\rm tr}(A^2P_{r})f\right]+\langle\left(\frac{\partial
X}{\partial t}\right)^{\top},\nabla S_{r+1}\rangle.
\end{equation}
\end{lemma}

\begin{proof}
Formula (\ref{eq:fundamental_relation_Reilly}) gives
$$(r+1)S_{r+1}=(-1)^r{\rm tr}(AP_r)=\displaystyle(-1)^r\sum_{i,j}a_{ji}a_{ij}^r =
\frac{1}{r!}\displaystyle\sum_{i,j,i_k,j_k}\epsilon^{j_1...j_rj}_{i_1...i_ri}a_{j_1i_1}...a_{j_ri_r}a_{ji},$$
where the functions $S_{r+1}$ are seen as functions of $t$. So,
differentiation with respect to $t$

\begin{eqnarray*}
(r+1)S'_{r+1}&=&\frac{1}{r!}\displaystyle\sum_{i,j,i_k,j_k}\epsilon^{j_1...j_kj}_{i_1...i_ki}[a'_{j_1i_1}...a_{j_ri_r}a_{ji}+a_{j_1i_1}...a_{j_ri_r}a'
_{ji}]\\
&=&\frac{(r+1)}{r!}\displaystyle\sum_{i,j,i_k,j_k}\epsilon^{j_1...j_rj}_{i_1...i_ki}a'_{ji}a_{j_1i_1}...a_{j_ri_r}\\
&=&(-1)^r\displaystyle\sum_{i,j}a'_{ji}a_{ij}^r=(-1)^r{\rm
tr}\left(\frac{\partial A}{\partial t}P_r\right).\end{eqnarray*}

Now, it is enough to compute $(-1)^r{\rm tr}\left(\frac{\partial
A}{\partial t}P_r\right)$.

\begin{eqnarray*}
S'_{r+1}&=&(-1)^r{\rm tr}\left(\frac{\partial A}{\partial t}P_r\right)=(-1)^r\sum_k\langle\frac{\partial A}{\partial t}P_re_k,e_k\rangle=\sum_k\langle 
S_r(A_k)\langle (\overline{\nabla}_{\frac{\partial X}{\partial t}}A)e_k,e_k\rangle\\
&=&\sum_k S_r(A_k)\left[\langle \overline{\nabla}_{\frac{\partial X}{\partial t}}Ae_k,e_k\rangle-\langle A\overline{\nabla}_{\frac{\partial 
X}{\partial t}}e_k,e_k\rangle\right]\\
&=&-\sum_k S_r(A_k)\langle \overline{\nabla}_{\frac{\partial
X}{\partial t}}\overline{\nabla}_{e_k}N,e_k\rangle- \sum_k
S_r(A_k)\langle A\overline{\nabla}_{e_k}\partial X/\partial
t,e_k\rangle,\end{eqnarray*} where we used that $[\partial
X/\partial t,e_k]=0$ in the last term.

Now, if $\overline{R}$ denotes the curvature tensor of
$\overline{M}$, we have
$$\overline{R}(e_k,\partial X/\partial t)N=\overline{\nabla}_{\frac{\partial X}{\partial 
t}}\overline{\nabla}_{e_k}N-\overline{\nabla}_{e_k}\overline{\nabla}_{\frac{\partial X}{\partial t}}N+\overline{\nabla}_{[e_k,\frac{\partial 
X}{\partial t}]}N.$$ Thus, by using also (\ref{eq:decomposition of
variational vector field})
\begin{eqnarray*}
S'_{r+1}
&=&-\sum_k S_r(A_k)\left[\langle\overline{R}(e_k,\partial X/\partial t)N,e_k\rangle+\langle\overline{\nabla}_{e_k}\overline{\nabla}_{\frac{\partial 
X}{\partial t}}N,e_k\rangle\right]\\
&&-\sum_kS_r(A_k)\langle\overline{\nabla}e_{k}(fN+(\partial
X/\partial t)^{\top}),Ae_k\rangle.\end{eqnarray*}

Since that the ambient spacetime is of constant sectional curvature,
it yields that
$$\langle\overline R(X,Y)W,Z\rangle=c\{\langle X,W\rangle\langle
Y,Z\rangle-\langle X,Z\rangle\langle Y,W\rangle\}.$$ Then
\begin{eqnarray*}
S'_{r+1}&=&-\sum_k S_r(A_k)c(\langle e_k,N\rangle\langle\partial X/\partial t,e_k\rangle-\langle e_k,e_k\rangle\langle\partial X/\partial t,N\rangle)\\
&&-\sum_k S_r(A_k)\langle\overline{\nabla}_{e_k}\overline{\nabla}_{\frac{\partial X}{\partial t}}N,e_k\rangle-\sum_k S_r(A_k)\langle 
Ae_k,\overline{\nabla}e_{k}fN\rangle\\
&&-\sum_k S_r(A_k)\langle\overline{\nabla}_{e_k}(\partial X/\partial t)^{\top},Ae_k\rangle\\
&=&-c\sum_k S_r(A_k)f-\sum_k S_r(A_k)e_k\langle\overline{\nabla}_{\frac{\partial X}{\partial t}}N,e_k\rangle+\sum_k 
S_r(A_k)\langle\overline{\nabla}_{\frac{\partial X}{\partial t}}N,\overline{\nabla}_{e_k}e_k\rangle\\
&&-\sum_k S_r(A_k)\langle Ae_k,f\overline{\nabla}e_{k}N\rangle-\sum_k S_r(A_k)e_k\langle Ae_k,(\partial X/\partial t)^{\top}\rangle\\
&&+\sum_k S_r(A_k)\langle\overline{\nabla}_{e_k} Ae_k,(\partial
X/\partial t)^{\top}\rangle.\end{eqnarray*}

Now, by using the expression for the trace of the operator $P_r$, we
get
\begin{eqnarray*}
S'_{r+1}&=&-(-1)^rc{\rm tr}(P_r)f+\sum_k (-1)^rP_re_k\langle N,\overline{\nabla}_{\frac{\partial X}{\partial t}}e_k\rangle\\
&&-\sum_k S_r(A_k)\langle\overline{\nabla}_{\frac{\partial
X}{\partial t}}N,N\rangle\langle N,\overline{\nabla}_{e_k}e_k\rangle
+f\sum_k S_r(A_k)\langle Ae_k,Ae_k\rangle\\&&-\sum_k (-1)^rP_re_k\langle Ae_k,(\partial X/\partial t)^{\top}\rangle+\sum_k S_r(A_k)\langle\nabla_{e_k} Ae_k,(\partial X/\partial t)^{\top}\rangle\\
&=&-(-1)^rc{\rm tr}(P_r)f+\sum_k (-1)^rP_re_k\langle N,\overline{\nabla}_{e_k}\partial X/\partial t\rangle+f\sum_k(-1)^r\langle AP_re_k,Ae_k\rangle\\
&&-\sum_k (-1)^rP_re_k\langle Ae_k,\partial X/\partial t\rangle+\sum_k (-1)^r\langle\nabla_{P_re_k} Ae_k,(\partial X/\partial t)^{\top}\rangle\\
&=&-(-1)^rc{\rm tr}(P_r)f-\sum_k(-1)^r\left(P_re_ke_k(f)+P_re_k\langle\overline{\nabla}_{e_k}N,\partial X/\partial t\rangle\right)\\
&&+(-1)^r{\rm tr}(A^2P_r)f-\sum_k (-1)^rP_re_k\langle Ae_k,\partial X/\partial t\rangle\\
&&+(-1)^r\langle\sum_k \nabla_{P_re_k} Ae_k,(\partial X/\partial
t)^{\top}\rangle.
\end{eqnarray*}

Now, by using Codazzi's equation, we have
\begin{eqnarray*}
\sum_k\nabla_{P_re_k}Ae_k&=&\sum_k\left(\nabla_{e_k}AP_re_k+A[P_re_k,e_k]\right)=\sum_k\nabla_{e_k}((-1)^rS_{r+1}I+P_{r+1})e_k\\
&&+\sum_k A\left(\nabla_{P_re_k}e_k-\nabla_{e_k}P_re_k\right)\\
&=&\sum_k(-1)^re_k(S_{r+1})e_k+\sum_k\nabla_{e_k}P_{r+1})e_k-A\nabla_{e_k}P_re_k\\
&=&(-1)^r{\rm div} S_{r+1}+{\rm div}(P_{r+1})-A({\rm div}P_r)=\nabla
S_{r+1},
\end{eqnarray*}
since that the operators $P_r$ are free of divergence. Hence
$$S'_{r+1}=(-1)^{r+1}\left(c{\rm tr}(P_r)f+L_{r}f-{\rm tr}(A^2P_r)f\right)+\langle\nabla S_{r+1},(\partial X/\partial t)^{\top}\rangle.$$

\end{proof}

The previous Lemma allows us to compute the first variation of the $r-$area functional.

\begin{proposition}\label{prop:first variation}
Under the hypotheses of Lemma~\ref{lemma:computing the time-derivative of Sr}, if $X$ is a variation of $x$, then
\begin{equation}\label{eq:first variation}
\mathcal A_r'(t)=\int_M[(-1)^{r+1}(r+1)S_{r+1}+c_r]f\,dM_t,
\end{equation}
where $c_r=0$ if $r$ is even and
$c_r=-\frac{n(n-2)(n-4)\ldots(n-r+1)}{(r-1)(r-3)\ldots
2}(-c)^{(r+1)/2}$ if $r$ is odd.
\end{proposition}

\begin{proof}
We make an inductive argument. The case $r=0$ is well known, and to the case $r=1$ we use the classical formula
$$\frac{\partial}{\partial_t}dM_t=[-S_1f+{\rm div}(\partial X/\partial t)^{\top}]dM_t$$
to get
\begin{eqnarray*}
\mathcal A_1'&=&\int_M F_1'dM_t+\int_M F_1\frac{\partial}{\partial_t}dM_t\\
&=&-\int_M S_1'dM_t-\int_M F_1[-S_1f+{\rm div}(\partial X/\partial t)^{\top}]dM_t\\
&=&\int_M[\Delta f-(S_1^2-2S_2)f+ncf-\langle(\partial X/\partial t)^{\top},\nabla S_1\rangle\\
&&+S_1^2f-S_1{\rm div}(\partial X/\partial t)^{\top}]dM_t\\
&=&\int_M 2S_2fdM_t+nc\int_MfdM_t-\int_M{\rm div}\left(S_1(\partial X/\partial t)^{\top}\right)dM_t\\
&=&\int_M(2S_2+nc)fdM_t,
\end{eqnarray*}
where in the last equality we used that $M$ is closed and $X$ is volume-preserving.

Now, if $r\geq 2$, the induction hypothesis and (\ref{eq:differentiation of $S_r$}) give
\begin{eqnarray*}
\mathcal A_r'&=&\int_M F_r'dM_t+\int_M F_r\frac{\partial}{\partial_t}dM_t\\
&=&\int_M\left[(-1)^rS_r'-\frac{c(n-r+1)}{r-1}F_{r-2}'\right]dM_t\\
&&+\int_M\left((-1)^rS_r-\frac{c(n-r+1)}{r-1}F_{r-2}\right)\frac{\partial}{\partial_t}dM_t\\
&=&\int_M(-1)^r\left\{S_r'+S_r[-S_1f+{\rm div}\left(\partial X/\partial t\right)^T]\right\}dM_t\\
&&-\frac{c(n-r+1)}{r-1}\left[\int_MF_{r-2}'dM_t+\int_MF_{r-2}\frac{\partial}{\partial_t}dM_t\right]\\
&=&\int_M(-1)^r\left\{S_r'-S_1S_rf+S_r{\rm div}\left(\partial X/\partial t\right)^T\right\}dM_t -\frac{c(n-r+1)}{r-1}A_{r-2}'\\
&=&\int_M\left[c{\rm tr}(P_{r-1})f+L_{r-1}f-{\rm tr}(A^2 P_{r-1})f+(-1)^r\langle\nabla S_r,(\partial X/\partial t)^T\rangle\right] dM_t\\
&&+(-1)^r\int_M\left(-S_1S_rf+S_r{\rm div}\left(\partial X/\partial t\right)^T\right)dM_t-\frac{c(n-r+1)}{r-1}A_{r-2}'\\
&=&\int_M\left[c(-1)^{r-1}(n-r+1)S_{r-1}f-(-1)^{r-1}(S_1S_r-(r+1)S_{r+1})f\right]dM_t\\
&&+(-1)^r\int_M\langle\nabla S_r,(\partial X/\partial t)^{\top}\rangle dM_t\\
&&+\int_M\left((-1)^{r+1}S_1S_rf+(-1)^rS_r{\rm div}\left(\partial X/\partial t\right)^{\top}\right)dM_t\\
&&-\frac{c(n-r+1)}{r-1}\int_M[(-1)^{r-1}(r-1)S_{r-1}+c_{r-2}]fdM_t\\
&=&\int_M[(-1)^{r+1}(r+1)S_{r+1}-\frac{c(n-r+1)}{r-1}c_{r-2}]fdM_t\\
&&+(-1)^r\int_M{\rm div}(S_r(\partial X/\partial t)^{\top})dM_t.
\end{eqnarray*}
It now suffices to apply the divergence theorem and note that
$c_r=-\frac{c(n-r+1)}{r-1}c_{r-2}$.
\end{proof}

In order to characterize spacelike immersions of constant $(r+1)-$th mean curvature, let $\lambda$ be a real constant and $\mathcal 
J_r:(-\epsilon,\epsilon)\rightarrow\mathbb R$ be the {\em Jacobi functional} associated to the variation $X$, i.e.,
$$\mathcal J_r(t)=\mathcal A_r(t)-\lambda\mathcal V(t).$$
As an immediate consequence of (\ref{eq:first variation}) we get
$$\mathcal J_r'(t)=\int_M[b_rH_{r+1}+c_r-\lambda]fdM_t,$$
where $b_r=(r+1){n\choose r+1}$. Therefore, if we choose $\lambda=c_r+b_r\overline H_{r+1}(0)$, where
$$\overline H_{r+1}(0)=\frac{1}{\mathcal A_0(0)}\int_MH_{r+1}(0)dM$$
is the mean of the $(r+1)-$th curvature $H_{r+1}(0)$ of $M$, we arrive at
$$\mathcal J_r'(t)=b_r\int_M[H_{r+1}-\overline H_{r+1}(0)]fdM_t.$$
Hence, a standard argument (cf.~\cite{BdC:84}) shows that $M$ is a critical point of $\mathcal J_r$ for all variations of $x$ if and only if $M$ has 
constant $(r+1)-$th mean curvature.

We wish to study spacelike immersions $x:M^n\rightarrow\overline M^{n+1}$ that maximize $\mathcal J_r$ for all variations $X$ of $x$. The above 
dicussion shows that $M$ must have constant $(r+1)-$th mean curvature and, for such an $M$, leads us naturally to compute the second variation of 
$\mathcal J_r$. This, in turn, motivates the following

\begin{definition}
Let $\overline M^{n+1}_c$ be a Lorentz manifold of constant
sectional curvature $c$, and $x:M^n\rightarrow\overline M^{n+1}$ be
a closed spacelike hypersurface having constant $(r+1)-$th mean
curvature. We say that $x$ is strongly $r$-stable if, for every
smooth function $f:M\rightarrow\mathbb R$ one has $\mathcal
J_r''(0)\leq 0$.
\end{definition}

The sought formula for the second variation of $\mathcal J_r$ is another straightforward consequence of Proposition~\ref{prop:first variation}.

\begin{proposition} Let $x:M^n\rightarrow\overline M^{n+1}_c$ be a closed spacelike hypersurface
of the time-oriented Lorentz manifold $\overline M^{n+1}_c$, having constant $(r+1)-$mean curvature $H_{r+1}$. If
$X:M^n\times(-\epsilon,\epsilon)\rightarrow\overline M^{n+1}_c$ is a variation of $x$, then
\begin{equation}\label{eq:second formula of variation}
\mathcal J_r''(0)=(r+1)\int_M\left[L_r(f)+c{\rm tr}(P_r)f-{\rm tr}(A^2P_r)f\right]fdM.
\end{equation}
\end{proposition}

\section{$r-$stable spacelike hypersurfaces in GRW's}

As in the previous section, let $\overline M^{n+1}$ be a Lorentz manifold.
A vector field $V$ on $\overline M^{n+1}$ is said to be {\em conformal} if
\begin{equation}
\mathcal L_V\langle\,\,,\,\,\rangle=2\psi\langle\,\,,\,\,\rangle
\end{equation}
for some function $\psi\in C^{\infty}(\overline M)$, where $\mathcal L$ stands for the Lie
derivative of the Lorentz metric of $\overline M$. The function $\psi$ is called the
{\em conformal factor} of $V$.

Since $\mathcal L_V(X)=[V,X]$ for all $X\in\mathcal X(\overline M)$, it follows from the
tensorial character of $\mathcal L_V$ that $V\in\mathcal X(\overline M)$ is conformal if
and only if
\begin{equation}\label{eq:1.1}
\langle\overline\nabla_XV,Y \rangle+\langle X,\overline\nabla_YV\rangle=2\psi\langle X,Y\rangle,
\end{equation}
for all $X,Y\in\mathcal X(\overline M)$. In particular, $V$ is a Killing vector field
relatively to $\overline g$ if and only if $\psi\equiv 0$.

Any Lorentz manifold $\overline M^{n+1}$, possessing a globally defined, timelike conformal
vector field is said to be a {\em conformally stationary spacetime}.

In what follows we need a formula first derived in~\cite{Alias:07}.
As stated below, it is the Lorentz version of the one stated and
proved in~\cite{Barros:09}.

\begin{lemma}\label{lemma:Lr of conformal vector field}
Let $\overline M^{n+1}_c$ be a conformally stationary Lorentz
manifold having constant sectional curvature $c$ and conformal
vector field $V$. Let also $x:M^n\rightarrow\overline M^{n+1}_c$ be
a spacelike hypersurface of $\overline M^{n+1}_c$ and $N$ be a
future-pointing, unit normal vector field globally defined on $M^n$.
If $\eta=\langle V,N\rangle$, then
\begin{eqnarray}\label{eq:Laplacian formula_I}
L_r\eta&=&{\rm tr}(A^2P_r)\eta-c\,{\rm tr}(P_r)\eta-b_rH_rN(\psi)\\
&&+b_rH_{r+1}\psi+\frac{b_r}{r+1}\langle V,\nabla
H_{r+1}\rangle,\nonumber
\end{eqnarray}
where $\psi:\overline M^{n+1}\rightarrow\mathbb R$ is the conformal
factor of $V$, $H_j$ is the $j-$th mean curvature of $x$ and $\nabla
H_j$ stands for the gradient of $H_j$ on $M$.
\end{lemma}

A particular class of conformally stationary spacetimes is that of
{\em generalized Robertson-Walker} spacetimes, or {\em GRW} for
short (cf.~\cite{ABC:03}), namely, warped products $\overline
M^{n+1}=I\times_{\phi}F^n$, where $I\subset\mathbb R$ is an interval
with the metric $-ds^2$, $F^n$ is an $n$-dimensional Riemannian
manifold and $\phi:I\rightarrow\mathbb R$ is positive and smooth.
For such a space, let $\pi_I:M^{n+1}\rightarrow I$ denote the
canonical projection onto the $I-$factor. Then the vector field
$$V=(\phi\circ\pi_I)\frac{\partial}{\partial s}$$
is a conformal, timelike and closed (in the sense that its dual
$1-$form is closed) one, with conformal factor $\psi=\phi'$, where
the prime denotes differentiation with respec to $s$. Moreover
(cf.~\cite{Montiel:99}), for $s_0\in I$, the (spacelike) leaf
$M_{s_0}^n=\{s_0\}\times F^n$ is totally umbilical, with umbilicity
factor $-\frac{\phi'(s_0)}{\phi(s_0)}$ with respect to the
future-pointing unit normal vector field $N$.

If $\overline M^{n+1}=I\times_{\phi}F^n$ is a GRW and
$x:M^n\rightarrow\overline M^{n+1}$ is a complete spacelike
hypersurface of $\overline M^{n+1}$, such that $\phi\circ\pi_I$ is
limited on $M$, then $\pi_F\big|_M:M^n\rightarrow F^n$ is
necessarily a covering map (cf.~\cite{ABC:03}). In particular, if
$M^n$ is closed then $F^n$ is automatically closed.

Also, recall (cf.~\cite{ONeill:83}) that a GRW as above has constant
sectional curvature $c$ if and only if $F$ has constant sectional
curvature $k$ and the warping function $\phi$ satisfies the ODE
$$\frac{\phi''}{\phi}=c=\frac{(\phi')^2+k}{\phi^2}.$$

We can now state and prove our main result, which generalizes the main theorems of~\cite{BBC:08} and~\cite{Liu:08}.

\begin{theorem}\label{thm:r_estability}
Let $\overline M^{n+1}_c=I\times_{\phi}F^n$ be a generalized
Robertson-Walker spacetime of constant sectional curvature $c$, and
$x:M^n\rightarrow\overline M^{n+1}_c$ be a closed strongly
$r-$stable spacelike hypersurface. If the warping function $\phi$ is
such that $H_r\phi''\geq\max\{H_{r+1}\phi',0\}$ and the set where
$\phi'=0$ has empty interior on $M$, then either $M^n$ is
$r-$maximal or a totally umbilical slice $\{s_0\}\times F$.
\end{theorem}

\begin{proof}
Since $M^n$ is strongly $r$-stable then
$$0\geq\mathcal J_r''(0)=(r+1)\int_M\left[L_r(f)+c{\rm tr}(P_r)f-{\rm tr}(A^2P_r)f\right]fdM$$
for all smooth $f:M\rightarrow\mathbb R$. In particular, if
$f=\eta$, where (as in Lemma~\ref{lemma:Lr of conformal vector
field}) $\eta=\langle V,N\rangle=\phi\langle\frac{\partial}{\partial
s},N\rangle$, and we take into account that $H_{r+1}$ is constant on
$M$, then
$$L_r\eta+c\,{\rm tr}(P_r)\eta-{\rm tr}(A^2P_r)\eta=-b_rH_rN(\phi')+b_rH_{r+1}\phi',$$
so that
\begin{equation}\label{eq:auxiliar para prova}
\int_M\left[-H_rN(\phi')+H_{r+1}\phi'\right]\phi\langle\frac{\partial}{\partial s},N\rangle dM\leq 0.
\end{equation}

Now, observe that
$\overline\nabla\phi'=-\langle\overline\nabla\phi',\frac{\partial}{\partial
s}\rangle\frac{\partial}{\partial s}=-\phi''\frac{\partial}{\partial
s}$, and hence
$$N(\phi')=\langle N,\overline\nabla\phi'\rangle=-\phi''\langle N,\frac{\partial}{\partial s}\rangle=\phi''\cosh\theta,$$
where $\theta$ is the hyperbolic angle between $N$ and
$\frac{\partial}{\partial s}$. Substituting the above into
(\ref{eq:auxiliar para prova}), we finally arrive at
$$\int_M\left[H_r\phi''\cosh\theta-H_{r+1}\phi'\right]\phi\cosh\theta dM\leq 0.$$
Arguing as in the end of the proof of Theorem 1.1 of~\cite{BBC:08}, we get
$$H_r\phi''(\cosh\theta-1)=0\ \ \text{and}\ \ H_r\phi''=H_{r+1}\phi'$$
on $M$. Since $H_{r+1}$ is constant on $M$, either $M$ is
$r-$maximal or $H_{r+1}\neq 0$ on $M$. If this last case happens,
the condition on the zero set of $\phi$ on $M$, together with the
above, gives $H_r\phi''\neq 0$ in a dense subset of $M$, and hence
$\cosh\theta=1$ on this set. By continuity, $\cosh\theta=1$ on $M$,
so that $M$ is a slice.
\end{proof}

The above result has an interesting application in the case in which
$\overline M^{n+1}_c$ is the de Sitter space of constant sectional
curvature $1$. For this, we make a brief description of this
spacetime. Let $\mathbb L^{n+2}$ denote the $(n+2)$-dimensional
Lorentz-Minkowski space ($n\geq 2$), that is, the real vector space
$\mathbb R^{n+2}$, endowed with the Lorentz metric
$$\left\langle v,w\right\rangle ={\displaystyle\sum\limits_{i=1}^{n+1}}v_{i}w_{i}-v_{n+2}w_{n+2},$$
for all $v,w\in\mathbb{R}^{n+2}$. We define the
$\left(n+1\right)$-dimensional de Sitter space $\mathbb S_1^{n+1}$
as the following hyperquadric of $\mathbb{L}^{n+2}$
$$\mathbb S_1^{n+1}=\left\{p\in L^{n+2};\left\langle p,p\right\rangle=1\right\}.$$
From the above definition it is easy to show that the metric induced
from $\left\langle\,\,,\,\right\rangle $ turns $\mathbb S_1^{n+1}$
into a Lorentz manifold with constant sectional curvature $1$.

Choose a unit timelike vector $a\in\mathbb L^{n+2}$, then
$V(p)=a-\langle p, a\rangle p$, $p\in\mathbb S_1^{n+1}$ is a
conformal and closed timelike vector field. It foliates the de
Sitter space by means of umbilical round spheres $M_{\tau}=\{p\in
\mathbb S_1^{n+1}:\langle p,a\rangle=\tau\}$, $\tau\in\mathbb R$.
The level set given by $\{p\in\mathbb S_1^{n+1}:\langle
p,a\rangle=0\}$ defines a round sphere of radius one which is a
totally geodesic hypersurface in $\mathbb S_1^{n+1}$. We will refer
to that sphere as the equator of $\mathbb S_1^{n+1}$ determined by
$a$.

In the context of warped products, the de Sitter space can be thought of as the following GRW
$$\mathbb S_1^{n+1}=-\mathbb R\times_{\cosh s}\mathbb{S}^n,$$
where $\mathbb S^n$ means Riemannian unit sphere. We observe that
there is a lot of possible choices for the unit timelike vector
$a\in\mathbb L^{n+2}$ and, hence, a lot of ways to describe $\mathbb
S_1^{n+1}$ as such a GRW (cf.~\cite{Montiel:99}, Section $4$). We
notice that in this model, the equator of $\mathbb S_1^{n+1}$ is the
slice $\{0\}\times\mathbb{S}^n$ and, consequently, $\phi'(s)=\sinh
s$ vanishes only in this slice.

In this setting, from Theorem~\ref{thm:r_estability}, we obtain the
following

\begin{corollary}
Let $x:M^n\rightarrow\mathbb S_1^{n+1}$ be a closed strongly
$r-$stable spacelike hypersurface. Suppose that the set of points in
which $M^n$ intersects the equator of $\mathbb S_1^{n+1}$ has empty
interior. If $H_r\geq\max\{H_{r+1},0\}$, then either $M^n$ is
$r-$maximal or a umbilical round sphere.
\end{corollary}

\end{document}